\documentclass[a4paper]{article}
\usepackage{amssymb,amsmath}
\usepackage{nth}
\usepackage{url}
\usepackage{hyperref}

\title{Semiarcs with a long secant in $\mathrm{PG}(2,q)$}
\author{Bence Csajb\'{o}k\thanks{Author was supported by the Hungarian National Foundation for Scientific Research, Grant No. K 81310.\newline$^\dag$ Author was supported by ERC Grant No. 227701 DISCRETECONT.\newline$^\ddag$ Author was supported by by the Slovenian--Hungarian Intergovernmental Scientific and 
Technological Cooperation Project, Grant No. T\'ET 10-1-2011-0606.}, \, Tam\'{a}s H\'{e}ger$^{*\dag}$ and Gy\"{o}rgy Kiss$^{*\ddag}$}
\date{}

\begin{document}
\maketitle

\newcommand{\cR}{{\mathcal R}}
\newcommand{\cC}{{\mathcal C}}
\newcommand{\cH}{{\mathcal H}}
\newcommand{\cN}{{\mathcal N}}
\newcommand{\cM}{{\mathcal M}}
\newcommand{\cE}{{\mathcal E}}
\newcommand{\cU}{{\mathcal U}}
\newcommand{\cK}{{\mathcal K}}
\newcommand{\cP}{{\mathcal P}}
\newcommand{\cL}{{\mathcal L}}
\newcommand{\cQ}{{\mathcal Q}}
\newcommand{\cX}{{\mathcal X}}
\newcommand{\F}{{\mathbf F}}
\newcommand{\GF}{\hbox{{\rm GF}}}
\newcommand{\st}{{\cal S}_t} 
\newcommand{\so}{{\cal S}} 
\newcommand{\sk}{{\cal S}_2} 
\newcommand{\ad}{E}
\newcommand{\lmb}{\lambda}
\newcommand{\ep}{\varepsilon}
\newcommand{\la}{\langle}
\newcommand{\ra}{\rangle}
\newcommand{\Qed}{\hfill $\Box$ \medskip}
\newcommand{\cB}{{\mathcal B}}
\newcommand{\cT}{{\mathcal T}}
\newcommand{\ind}{\mathrm{ind}}
\newcommand{\PG}{\mathrm{PG}}
\newcommand{\AG}{\mathrm{AG}}
\newcommand{\os}{\so}

\newtheorem{theorem}{Theorem}[section]
\newtheorem{corollary}[theorem]{Corollary}
\newtheorem{lemma}[theorem]{Lemma}
\newtheorem{example}[theorem]{Example}
\newtheorem{proposition}[theorem]{Proposition}
\newtheorem{remark}[theorem]{Remark}
\newtheorem{definition}{Definition}

\newcommand{\R}{{\mathbb R}}
\newcommand{\N}{{\mathbb N}}
\newcommand{\Z}{{\mathbb Z}}
\newcommand{\proof}{\noindent{\bf Proof.}\ }

\begin{abstract}
A $t$-semiarc is a pointset ${\cal S}_t$ with the property that the number of tangent lines to ${\cal S}_t$ at each
of its points is $t$. We show that if a small $t$-semiarc ${\cal S}_t$ in $\mathrm{PG}(2,q)$ has a large collinear subset ${\cal K}$, then
the tangents to ${\cal S}_t$ at the points of ${\cal K}$ can be blocked by $t$ points not in ${\cal K}$. 
We also show that small $t$-semiarcs are related to certain small blocking sets. 
Some characterization theorems for small semiarcs  with large collinear subsets in $\mathrm{PG}(2,q)$ are given.
\end{abstract}



\section{Introduction}
\label{sec:Intro}
\indent

Ovals, $k$-arcs and semiovals of finite projective planes are not only 
interesting geometric structures,
but they have applications to coding theory and cryptography, too \cite{bat}.
For details about these objects we refer the reader to 
\cite{jwph1, kgy2}.   

Semiarcs are natural generalizations of arcs. Throughout the paper $\Pi_q$ denotes an arbitrary projective
plane of order $q$.
A non-empty pointset ${\cal S}_t\subset \Pi _q$ is called a {\it t-semiarc} 
if for every point $P\in {\cal S}_t$
there exist exactly $t$ lines $\ell _1,\ell _2,\ldots \ell _t$ such that 
${\cal S}_t \cap \ell_i = \{P\}$ for $i=1,2,\ldots ,t$. 
These lines are called the \emph{tangents} to ${\cal S}_t$
at $P$. If a line $\ell $ meets $\st $ in $k$ points, then $\ell$ is called
a \emph{$k$-secant} of $\st$. We say that a $k$-secant is \emph{long}, if $q-k$ is a small number 
(which will be given a precise meaning later). 
The classical examples of $t$-semiarcs are the $k$-arcs (with $t=q+2-k$), 
subplanes (with $t=q-m$, where $m$ is the order of the subplane) and 
semiovals (that is semiarcs with $t=1$).

Because of the huge diversity of semiarcs, the complete classification 
is hopeless. The aim of this paper is to investigate and characterize 
semiarcs having some additional properties.
In Section~\ref{sec:direct} we consider a very special class,
namely $t$-semiarcs of size $k+q-t$ having a $k$-secant. These pointsets
are closely related to the widely studied structures defining few directions
\cite{Ball, BBSSZ, DirSzonyi}.
In Section~\ref{sec:main} 
we prove that in $\mathrm{PG}(2,q)$ if a small $t$-semiarc has a large collinear 
subset $\cK$, then the tangent 
lines at the points of $\cK$ belong to $t$ 
pencils, whose carriers are not in $\cK$.
This result generalizes the main result in Kiss \cite{kgy}. 
Small semiovals with large collinear subsets were studied in arbitrary projective planes as well, see 
Bartoli \cite{bar} and Dover \cite{dov}. 
The essential part of our proof is algebraic, it is
based on an application of the R\'edei polynomial and the 
Sz\H{o}nyi--Weiner Lemma.
In Section~\ref{sec:blsets} we associate to each $t$-semiarc
$\st$ a blocking set. If $\st$ is small and has a long secant, then
the associated blocking set is small. Applying theorems about the structure
of small blocking sets we prove some characterization theorems for
semiarcs. 

When $t\geq q-2$, then it is easy to characterize $t$-semiarcs.
If $t=q+1,\,q$ or $q-1$, then $\st$ is single point, a subset of a line,
or three non-collinear points, respectively; see \cite[Proposition 2.1]{csk}. 
Hence, if no other bound is specified, we usually will assume that $t\leq q-2$. 
If $t=q-2$, then it follows from \cite[Proposition 3.1]{csk} that 
$\st$ is one of the following three configurations:  
four points in general position, the six vertices of a complete 
quadrilateral, or a Fano subplane. Thus sometimes we may assume that $t\leq q-3$, which we indicate individually. 

Throughout the paper we use the following notation.
We denote points at infinity of $\PG(2,q)$, i.e. points on the line $\ell_{\infty}=[0:0:1]$, by $(m)$
instead of the homogeneous coordinates $(1:m:0)$.
We simply write $Y_{\infty}$ and $X_{\infty}$ instead of $(0:1:0)$ and $(1:0:0)$, respectively.
The points of $\ell_{\infty}$ are also called directions.
For affine points, i.e. points of $\PG(2,q)\setminus \ell_{\infty}$,  
we use the Cartesian coordinates $(a,b)$ instead of $(a:b:1)$.
If $P$ and $Q$ are distinct points in $\Pi_q$, then $PQ$ denotes the unique 
line joining them. If $\cal A$ and $\cal B$ are two pointsets in $\Pi_q$, 
then ${\cal A} \triangle {\cal B}$ denotes their symmetric difference, 
that is $({\cal A}\setminus{\cal B})\cup({\cal B}\setminus{\cal A})$.

Blocking sets play an important role in our proofs. For the sake of 
completeness we collect the basic definitions and some results about these objects.
A \emph{blocking set} $\cal B$ in a projective or affine plane is a set of points which intersects every line.
If $\cal B$ contains a line, then it is called trivial.
A point $P$ in a blocking set $\cal B$ is \emph{essential} 
if ${\cal B}\setminus \{P\}$ is not a blocking set, i.e.\ 
there is a tangent line to $\cal B$ at the point $P$.
A blocking set is said to be \emph{minimal} when no proper subset of it is 
a blocking set or, equivalently, each of its points
is essential. If $\ell$ is a line containing at most $q$ points of a blocking set $\cal B$ in $\Pi_q$, then
$|{\cal B}|\geq q+|\ell \cap \cB|$. In case of equality $\cal B$ is a blocking
set of \emph{R\'{e}dei type} and $\ell$ is a \emph{R\'{e}dei line} of $\cB$. 
Note that we also consider a line to be a blocking set of R\'{e}dei type. 
A blocking set in $\Pi _q$ is said to be \emph{small} if its size is less 
than $3(q+1)/2$.

\begin{theorem}[{\cite[Remark 3.3 and Corollary 4.8]{Szonyi}}]
\label{modp}
Let $\cB$ be a blocking set in $\mathrm{PG}(2,q)$, $q=p^h$, $p$ prime. 
If $|\cB|\leq 2q$, then $\cB$ contains a unique minimal blocking set.
If $\cB$ is a minimal blocking set of size less than $3(q+1)/2$, 
then each line intersects $\cB$ in $1 \pmod p$ points.
\end{theorem}

Note that a blocking set contains a unique minimal blocking set if and only 
if the set of its essential points is a blocking set.

\begin{theorem}[{\cite[Corollary 5.1]{linconj}, \cite{Polverino}, \cite{Szonyi}}]
\label{modpe}
Let $\cB$ be a minimal blocking set in $\mathrm{PG}(2,q)$, $q=p^h$, $p$ prime, 
of size $|\cB|<3(q+1)/2$. Then there exists a positive integer $e$, called the 
exponent of $\cB$, such that $e$ divides $h$, and
\[ q+1+p^e\left\lceil \frac{q/p^e+1}{p^e+1}\right\rceil
\leq |\cB| \leq \frac{1+(p^e+1)(q+1)-\sqrt{D}}{2}, \]
where $D=(1+(p^e+1)(q+1))^2-4(p^e+1)(q^2+q+1)$.

If $p^e\neq 4$ and $|\cB|$ lies in the interval belonging to $e$, 
then each line intersects $\cB$ in $1 \pmod {p^e}$ points.
\end{theorem}

\begin{theorem}[\cite{Blokhuis, smallb, bruen}]
\label{smallbs}
Let $\cB$ be a minimal blocking set in $\PG(2,q)$, $q=p^h$, $p$ prime.
Let $|\cB|=q+1+k$, and let $c_p=2^{-1/3}$ for $p=2,3$ and $c_p=1$ for $p>3$. Then the following hold.
\begin{enumerate}
	\item If $h=1$ and $k\leq(q+1)/2$, then $\cB$ is a line, or $k=(q+1)/2$ and each point of 
	$\cB$ has exactly $(q-1)/2$ tangent lines.
	\item If $h=2d+1$ and $k<c_pq^{2/3}$, then $\cB$ is a line.
	\item If $k\leq \sqrt{q}$, then $\cB$ is a line, or $k=\sqrt{q}$ and $\cB$ is a Baer subplane (that is a subplane of order $\sqrt{q}$). 
\end{enumerate}
\end{theorem}

We remark that the third point of the above theorem holds in arbitrary finite projective planes.

\section{Semiarcs and the direction problem}
\label{sec:direct}

If a $t$-semiarc $\st$ has a $k$-secant $\ell$, then its size $s$ is at least $k+q-t$,
because for any point $P\in \st \cap \ell $ there are $q+1-t$ non-tangent lines to $\st$
through $P$, one of which is $\ell$, and each of the remaining $q-t$ non-tangent lines contains 
at least one point from $\st \setminus \ell$. 
Thus we may always assume that $s=k+q-t+\ep$, where $\ep\geq 0$.
In this section we investigate the case $\ep =0$. Notice that $t<q$ implies $k\leq q+1-t$.

\begin{theorem}[{\cite[Theorem 4]{csb3}}]
\label{q+1-t}
In $\PG(2,q)$, a $t$-semiarc with a $(q+1-t)$-secant exists if and only if $t\geq (q-1)/2$.
\end{theorem}

\begin{proposition}[{\cite[Proposition 2.2]{csk}}]
\label{twolines}
Let $\Pi_q$ be a projective plane of order $q$, and let $t\leq q-2$. If a $t$-semiarc $\st$ in $\Pi_q$ is contained in the union of two lines, $\ell$ and $\ell'$, then 
$\ell\cap\ell'\notin\st$ and $|\ell\cap\st|=|\ell'\cap\st|=q-t$.
\end{proposition}

It is easy to give a combinatorial characterization of $t$-semiarcs of size $2(q-t)$ with a $(q-t)$-secant; 
for semiovals it was also proved by Bartoli \cite[Corollary 9]{bar}. 

\begin{proposition}
\label{komb}
Let $\Pi_q$ be a projective plane of order $q$, and let $t\leq q-2$.
If $\st$ is a $t$-semiarc of size $2(q-t)$ with a $(q-t)$-secant $\ell$, then
$\st$ consists of the symmetric difference of two lines with $t$ further points removed from each line.
\end{proposition}
\proof
Let $\cR=\st \setminus \ell$.
If $\ell'$ is a line joining two points of $\cR$,
then $\ell \cap \ell' \notin  \st$, otherwise there would be at least $t+1$ tangents to $\st$ at $\ell\cap\ell'$.
Now suppose to the contrary that there exist three non-collinear points in $\cR$. 
They determine three lines, each of which intersects $\ell$ in
$\ell \setminus \st$; hence at these three points of $\cR$ there are at most $t-1$ tangents to $\st$, a contradiction.
Thus the points of $\cR$ are collinear and $\ell \cap \ell' \notin \st$.
\Qed

The following example shows the existence of $t$-semiarcs of size
$k+q-t$ with three $k$-secants for odd values of $t$.

\begin{example}
\label{projtriangle0}
Let $C$ denote the set of non-squares in the field $\GF(q)$, $q$ odd.
The pointset $\{(0:1:s), (s:0:1), (1:s:0) \, : \, -s\in C\}$ is a semioval in $\PG(2,q)$ of size $3(q-1)/2$ 
with three $(q-1)/2$-secants (see Blokhuis \cite{seminuclear}).
If we delete $r<(q-1)/2-1$ points from each of the $(q-1)/2$-secants, then the remaining pointset is a $t$-semiarc 
of size $k+q-t$ with three $k$-secants, where $k=(q-1)/2-r$ and $t=2r+1$.
\end{example}

There also exist examples if $t$ is even. 
To give their construction, we need some notation.
A \emph{$(k,n)$-arc} is a set of $k$ points such that each line 
contains at most $n$ of these points.
A set $\cal T$ of $q+t$ points in $\Pi_q$ for which each line meets $\cal T$ in 0, 2 or $t$ points ($t\neq 0,2$) 
is either an oval (for $t=1$), or a $(q+t,\,t)$-arc of type $(0,\,2,\,t)$. 
Korchm\'{a}ros and Mazzocca \cite[Proposition 2.1]{KMarcs} proved that
$(q+t,\,t)$-arcs of type $(0,\,2,\,t)$ exist in $\Pi_q$ only if $q$ is even and $t\mid q$.
They also provided infinite families of examples in $\PG(2,q)$ whenever the field
$\GF(q/t)$ is a subfield of $\GF(q)$.
It is easy to see that
through each point of $\cal T$ there passes exactly one $t$-secant.
In \cite{t-arc} new constructions were given by G\'{a}cs and Weiner, and 
they proved that in $\PG(2,q)$ the 
$q/t+1$ $t$-secants of $\cal T$ pass through one point, called the $t$-nucleus of $\cal T$
(for $t=1$ and arbitrary projective plane of even order, see  \cite[Lemma 8.6]{jwph1}).
Recently Vandendriessche \cite{PV} found a new infinite family with $t=q/4$. 

\begin{example}
\label{redei2}
Let $\cal T$ be a $(q+\tau,\,\tau)$-arc of type $(0,\,2,\,\tau)$ in $\PG(2,q)$.
Delete $r<\tau-1$ points from each of the $\tau$-secants of $\cal T$. The remaining $k+q-t$ points form a $t$-semiarc with
$q/\tau+1$ $k$-secants, where $k=\tau-r$ and $t=rq/\tau$.
\end{example}

Since $(q+q/2,\,q/2)$-arcs of type $(0,\,2,\,q/2)$ exist, this construction gives $t$-semiarcs in $\PG(2,q)$, $q$ even, 
for each $t\leq q-4$, $t$ even. The following example is based on the combinatorial properties of subplanes.

\begin{example}
\label{subplane}
Let $\Pi_{\sqrt{q}}$ be a Baer subplane in the projective plane $\Pi_q$,
$q\geq 9$, and let $\ell$ be an extended line of $\Pi_{\sqrt{q}}$. Let $\cP$ 
be a set of $t\leq q-\sqrt{q}-2$ points in $\Pi_{\sqrt{q}}\setminus \ell$ such that no line intersects $\cP$ in exactly $\sqrt{q}-1$ points.
For example a $(t,\sqrt{q}-2)$-arc is a good choice for $\cP$.
Let $\cal T$ be a set of $t$ points in $\ell \setminus \Pi_{\sqrt{q}}$.
Then the pointset $\st:=(\Pi_{\sqrt{q}} \triangle \ell) \setminus ({\cal T} \cup \cP)$ 
is a $t$-semiarc of size $k+q-t$ with a $k$-secant, where $k=q-\sqrt{q}-t$.
\end{example}
\proof
Recall that a Baer subplane is a blocking set. Let $P\in\st$. If $P\in\ell$, then a line through $P$ is tangent to $\st$ if and only if it intersects $\Pi_{\sqrt{q}}$ in a point of $\cP$. If $P\in\Pi_{\sqrt{q}}$, then a line of $\Pi_{\sqrt{q}}$ through $P$ intersects $\st$ in at least $\sqrt{q}-(\sqrt{q}-2)=2$ points, and any other line through $P$ is tangent to $\st$ if and only if it intersects $\ell$ in $\cal{T}$. Thus there are exactly $t$ tangents to $\st$ at $P$.
\Qed

The so-called direction problem is closely related to $t$-semiarcs of size $k+q-t$ having a $k$-secant. 
We briefly collect the basic definitions and some results about this problem.
Consider $\mathrm{PG}(2,q)=\mathrm{AG}(2,q)\cup \ell _{\infty}$.
Let $\cU$ be a set of points of $\AG(2,q)$. 
A point $P$ of $\ell_{\infty}$ is called a direction determined by $\cU$ if 
there is a line through $P$ that contains at least two points of $\cU$. 
The set of directions determined by $\cU$ is denoted by $D_{\cU}$.
If $|\cU|=q$ and $Y_{\infty}\notin D_{\cU}$, then $\cU$ can be considered as a graph of a function, 
and $\cU \cup D_{\cU}$ is a blocking set of R\'{e}dei type.
Our next construction is based on the following result of Blokhuis et al.\ \cite{BBSSZ} and Ball \cite{Ball}.

\begin{theorem}[\cite{BBSSZ,Ball}]
\label{bbbssz}
Let $\cU\subset \mathrm{AG}(2,q)$, $q=p^h$, $p$ prime, be a pointset of size $q$. 
Let $z=p^e$ be maximal having the property that if $P\in D_{\cU}$ and $\ell$ is a line through 
$P$, then $\ell $ intersects $\cU$ in $0 \pmod z$ points. 
Then one of the following holds:
\begin{enumerate}
	\item $z=1$ and $(q+3)/2 \leq |D_{\cU}| \leq q+1$,
	\item $\GF(z)$ is a subfield of $\GF(q)$ and $q/z + 1 \leq |D_{\cU}| \leq (q - 1)/(z - 1)$,
	\item $z=q$ and $|D_{\cU}|=1$.
\end{enumerate}
\end{theorem}

Let $\cal B$ be a small blocking set of R\'{e}dei type in $\PG(2,q)$, $q=p^h$, $p$ prime, 
and let $\ell$ be one of its R\'{e}dei lines.
Since $|{\cal B}|<3(q+1)/2$, we have $|\ell \cap {\cal B}|<(q+3)/2$. 
Hence the previous theorem implies that there exists an integer $e$ such that 
$e$ divides $h$, $1<p^e\leq q$ holds and each affine line intersects $\cal B$ in $1 \pmod{p^e}$ points.
Starting from $\cal B$, we give a generalization of Example \ref{subplane}, 
which is also a semiarc for similar reasons. 

\begin{example}
\label{redei}
Let $\cal B$ be a small blocking set of R\'{e}dei type in $\PG(2,q)$ 
and let $\ell$ be one of its R\'{e}dei lines. Denote by $z=p^e$ the maximal number such that 
each line intersects $\cal B$ in $1 \pmod z$ points and suppose $z\geq 3$. 
Let $\cP$ be a set of $t\leq q-|{\cal B}\cap\ell|-1$ points in ${\cal B}\setminus \ell$ such that
for each line $\ell'$ intersecting $\cB$ in more than one point we 
have $|\ell' \cap \cP| \neq |\ell' \cap \cB|-2$.
For example a $(t,z-2)$-arc is a good choice for $\cP$.
Also let $\cal T$ be a set of $t$ points in $\ell \setminus {\cal B}$. Then 
the pointset $\st:=({\cal B} \triangle \ell) \setminus ({\cal T} \cup \cP)$ 
is a $t$-semiarc of size $k+q-t$ with a $k$-secant, where $k=2q+1-|{\cal B}|-t$.
\end{example}

Note that if $\cal B$ is a line, then Example \ref{redei} gives the example seen in Proposition \ref{komb}. 
To characterize the examples above, we need results about the number of 
directions determined by a set of $q$ affine points, and results about the 
extendability of a set of almost $q$ affine points to a set of $q$ points 
such that the two pointsets determine the same directions. The first theorem 
about the extendability was proved by Blokhuis \cite{seminuclear}; 
see also Sz\H{o}nyi \cite{DirSzonyi}.

\begin{theorem}[{\cite[Proposition 2]{seminuclear}, \cite[Remark 7]{DirSzonyi}}]
\label{Dir1}
Let $\cU \subset \AG(2,q)$, $q\geq 3$, be a pointset of size $q-1$. 
Then there exists a unique point $P$ such that the pointset $\cU \cup \{P\}$
determines the same directions as $\cU$.
\end{theorem}

\noindent
Slightly extending a result of Sz\H{o}nyi \cite[Theorem 4]{DirSzonyi}, Sziklai proved the following.

\begin{theorem}[{\cite[Theorem 3.1]{DirSziklai}}]
\label{Dir}
Let $\cU \subset \AG(2,q)$ be a pointset of size $q-n$, where $n\leq\alpha \sqrt{q}$ for some $1/2 \leq \alpha < 1$.
If $|D_{\cU}|<(q+1)(1-\alpha)$, then $\cU$ can be extended to a set $\cU'$ of size $q$ 
such that $\cU'$ determines the same directions as $\cU$.
\end{theorem}

\noindent
The three cases of the next theorem were proved by Lov\'{a}sz and 
Schrijver \cite{LS}, by G\'{a}cs \cite{Gacs}, and by G\'{a}cs, Lov\'{a}sz 
and Sz\H{o}nyi \cite{p2}, respectively.

\begin{theorem}[{\cite{LS,Gacs,p2}}]
\label{pp2}
Let $\cU $ be the set of $q$ affine points in $\AG(2,q)$, $q=p^h$, $p$ prime.
\begin{enumerate}
	\item If $h=1$ and $|D_{\cU}|=(p+3)/2$, then $\cU $ is affinely equivalent to the graph of the
function $x\mapsto x^{\frac{p+1}{2}}$.
	\item If $h=1$ and $|D_{\cU}|>(p+3)/2$, then $|D_{\cU}|\geq\lfloor 2(p-1)/3\rfloor+1$.
	\item If $h=2$ and $|D_{\cU}|\geq (p^2+3)/2$, then either 
$|D_{\cU}|=(p^2+3)/2$ and $\cU $ is affinely equivalent to the graph of the
function $x\mapsto x^{\frac{p^2+1}{2}}$, or $|D_{\cU}|\geq (p^2+p)/2+1$.
\end{enumerate}
\end{theorem}

For the characterization of semiarcs in Example \ref{redei} we also need the following lemma.

\begin{lemma}
\label{small0}
Let $z$ and $t$ be two positive integers such that $z\geq 3$ and $t \leq \sqrt{q(z-1)/z}$.
Also let $\cU \subset \AG(2,q)$ be a set of $q-t$ affine points
and let $E\subseteq F$ be two sets of directions satisfying the following properties:
\begin{enumerate}
	\item  
	There are  at least $t$ tangents to $\cU$ with direction in $F$ through each point of $\cU$;
	\item 
	there exists a suitable set of $t$ affine points, $\cP$, such that $\cU \cap \cP =\emptyset$ and each tangent to $\cU$  with direction not in 
	$E$ intersects $\cU \cup \cP$ in $0 \pmod z$ points.
\end{enumerate}
Then $|E|\geq t$.
\end{lemma}
\proof
If $\ell $ is a tangent to $\cU$ that intersects $F\setminus E$, then $|\cP \cap \ell|\equiv -1$ (mod $z$).
The maximum number of such tangent lines is $\frac{t(t-1)}{(z-1)(z-2)}$. Hence at least
$(q-t)t - \frac{t(t-1)}{(z-1)(z-2)}$ tangents to $\cU$ have direction in $E$.
This implies $$|E|q  \geq (q-t)t - \frac{t(t-1)}{(z-1)(z-2)}\mbox{, \quad thus \quad} 
(|E|-t)q \geq -t^2-\frac{t(t-1)}{(z-1)(z-2)}.$$
If $|E|-t$ is a negative integer, then this inequality gives
$q<t^2\frac{(z-1)(z-2)+1}{(z-1)(z-2)}\leq t^2z/(z-1)$, contradicting the assumption $t\leq \sqrt{q(z-1)/z}$.
\Qed

\begin{theorem}
\label{invsemiarc}
Let $\st$ be a $t$-semiarc in $\PG(2,q)$, $q=p^h$, $p$ prime, of size 
$k+q-t$ and let $\ell$ be a $k$-secant of $\st$. Then the conditions
\begin{itemize}
\item 
$t=1$, $q>4$ and $k > (q-1)/2$, or 
\item
$2\leq t \leq \alpha \sqrt{q} $ and $k>\alpha(q+1)$ for some $1/2 \leq \alpha \leq \sqrt{(p-1)/p}$ if $p$ is an odd prime, and
$1/2 \leq \alpha \leq \sqrt{3}/2$ if $p=2$
\end{itemize}
imply that $\st$ is a semiarc described in Example \ref{redei}.
\end{theorem}
\proof
Take $\ell$ as the line at infinity and let $\cU=\st \setminus \ell \subseteq \AG(2,q)$. 
The directions in $\st \cap \ell$ are not determined by $\cU$, 
hence $|D_{\cU}| < (q+1)(1-\alpha)$ holds for $t\geq 2$.
We can apply Theorem \ref{Dir1} when $t=1$; if $t\geq2$, then the conditions of 
Theorem \ref{Dir} hold since $|\cU|=q-t$ and $t\leq \alpha \sqrt{q}$.
Let $\cP=\{P_1,P_2,\ldots,P_t\}$ be a set of $t$ points such that $\cU\cup \cP$ determines 
the same directions as $\cU$.

First consider the case $t\geq 2$. We have $|D_{\cU}|<(q+1)/2$, thus applying Theorem \ref{bbbssz} 
we get that there exists an integer $z=p^e\geq 3$ such that each affine line with direction in $D_{\cU}$ intersects 
$\cU\cup \cP$ in $0 \pmod z$ points. 
We can apply Lemma \ref{small0} with $F=\ell \setminus \st$ and $E=\ell \setminus (\st \cup D_{\cU})$ to obtain
$|E|\geq t$. On the other hand the lines joining any point of $E$ with any point of $\cU$ are tangents to $\st$, 
thus $|E|\leq t$.
The same observation implies that each of the tangents to $\st$ at the points of $\cU$ meets $E$.
Let $\cB=\cU \cup \cP \cup D_{\cU}$, which is a small blocking set of R\'{e}dei type.
Let $\ell'\neq \ell$ be a line intersecting $\cB$ in more than one point and let $M=\ell' \cap \ell$.
Then $M \in D_{\cU}\subseteq \cB$ and $M\notin E$.
If $|\ell' \cap \cP|=|\ell' \cap \cB|-2$, then $\ell'$ would be a tangent to $\st$ at the unique point of $\ell' \cap \cU$, but this
is a contradiction since $M\notin E$.
We obtained Example \ref{redei}.

If $t=1$, then in the same way we get that there exists an integer $z=p^e\geq
2$ such that each affine line with direction in $D_{\cU}$ intersects 
$\cU\cup \{ P_1\} $ in $0 \pmod z$ 
points. 
If $z\geq 3$, then we can finish the proof as above, otherwise Theorem \ref{bbbssz} implies
$|D_{\cU}|\geq q/2+1$. Compared to $|D_{\cU}|<(q+3)/2$, we get $|D_{\cU}|=q/2+1$ and hence $k=q/2$. This means that each of the $q-1$ tangent lines at the points of $\cU$ passes through $P_1$. If $q>4$, then $q-1>q/2+1$, thus at least one of these tangents 
would intersect $\ell$ in $\st$, that is a contradiction.
\Qed

Next, as a corollary of Theorem \ref{pp2}, we get the characterization of the semioval $(t=1)$ cases of Examples \ref{projtriangle0} and \ref{subplane} in planes of prime or prime square order. 

\begin{corollary}
\label{cort1}
Let $\so_1$ be a semioval of size $k+q-1$ in $\PG(2,q)$, $3 \leq q=p^h$, $p$ prime,
$h\leq 2$, and let $\ell$ be a $k$-secant of $\so_1$.
Then we have the following.
\begin{enumerate}
\item
If $h=1$ and $k>(p+4)/3$, then there are two possibilities:
\begin{itemize}
\item
$k=q-1$ and $\so_1$ is the semioval characterized in Proposition \ref{komb},
\item 
$\so_1$ is the 
semioval described in Example \ref{projtriangle0}.
\end{itemize}

\item
If $h=2$ and $k>(p^2-p)/2$, then there are four possibilities:
\begin{itemize}
\item
$k=q-1$ and $\so_1$ is the semioval characterized in Proposition \ref{komb},
\item
$\so_1$ is the 
semioval described in Example \ref{projtriangle0}, 
\item
$\so_1$ is the 
semioval described in Example \ref{subplane}, 
\item
$p=2$ and $\so_1$ is an oval in $\mathrm{PG}(2,4)$.
\end{itemize}
\end{enumerate}
\end{corollary}
\proof
Consider $\ell$ as the line at infinity and let $\cU=\so_1 \setminus \ell$. 
The points of $\ell \cap \so_1$ are not determined directions, 
hence we have $k+|D_{\cU}|\leq q+1$.
As the pointset $\cU $ has size $q-1$, it follows
from Theorem \ref{Dir1} that there exists a point $P$ such that $\cU\cup \{P\}$ 
determines the same directions as $\cU$.

First consider the case $h=1$. If $k>(p+4)/3$, then 
$|D_{\cU}|<\lfloor 2(p-1)/3\rfloor +1$ and thus Theorems \ref{bbbssz} and \ref{pp2} imply that $|D_{\cU}|=1$ and $\cU$ is contained in a line,
or $|D_{\cU}|=(p+3)/2$ and $\cU \cup \{P\}$ is affinely equivalent to the 
graph of the function $x\mapsto x^{\frac{p+1}{2}}$. In the first case it is 
easy to see that $\so_1$ is the semioval characterized in Proposition
\ref{komb}. In the latter case the graph of $x\mapsto x^{\frac{p+1}{2}}$ 
is contained in two lines, namely $[1:1:0]$ and $[1:-1:0]$, 
and these lines are $(q-1)/2$-secants. Hence $\so_1$ is contained in the 
union of three lines and it has two $(q-1)/2$-secants. These semiovals were 
characterized by Kiss and Ruff \cite[Theorem 3.3]{smallsemi}; they proved that the 
only possibility is the semioval described in Example \ref{projtriangle0}.

Now suppose that $h=2$. If $k>(p^2-p)/2$, then $|D_{\cU}|<(p^2+p)/2+1$, thus $|D_{\cU}|\in \{1,(p^2+3)/2\}$ or $1<|D_{\cU}|<(p^2+3)/2$.
If $|D_{\cU}|=1$ or $|D_{\cU}|=(p^2+3)/2$, then we can argue as before.
In the remaining case it follows from
Theorems \ref{bbbssz} and \ref{smallbs} (or already from \cite[Theorem 5.7]{Szonyi}), 
that $|D_{\cU}|=p+1$ and $\cU\cup \{P\}\cup D_{\cU}$ 
is a Baer subplane.
If $p>2$, then $\so_1$ has exactly $p^2-p-k$ tangents at the points of 
$\cU$, hence $k=p^2-p-1$ and $\so_1$ is the semioval described in 
Example \ref{subplane}. Finally, if $p=2$, then $k\geq 2$ and $|D_{\cU}|=p+1=3$, thus $k=2$ and $\so$ is an oval in $\PG(2,4)$.
\Qed

\section{Proof of the main lemma}
\label{sec:main}

First we collect the most important properties of the 
{\it R\'edei polynomial.} Consider a subset 
$\cU=\{ (a_i,b_i):i=1,2,\ldots ,|\cU|\}$ of the affine plane $\AG(2,q)$.
The R\'edei polynomial of $\cU$ is 
\[H(X,Y)=\prod_{i=1}^{|\cU|}(X+a_iY-b_i)=\sum_{j=0}^{|\cU|}h_j(Y)X^{|\cU|-j}\in \GF(q)[X,Y],\]
where $h_j(Y)$ is a polynomial of degree at most $j$ in $Y$ and
$h_0(Y)\equiv 1$. Let $H_m(X)$ be the one-variable polynomial
$H(X,m)$ for any fixed value $m$. Then $H_m(X)\in {\rm GF}(q)[X]$ is 
a fully reducible polynomial, which reflects some geometric properties 
of $\cU$. 

\begin{lemma}[Folklore] Let $H(X,Y)$ be the R\'edei polynomial of the 
pointset $\cU$, and let $m\in {\rm GF}(q)$. Then $X=k$ is a root of
$H_m(X)$ with multiplicity $r$ if and only if the line with equation 
$Y=mX+k$ meets $\cU$ in exactly $r$ points.
\end{lemma} 

We need another result about polynomials. For $r\in \R,$ let $r^+=\max\{0,r\}$.

\begin{theorem}[Sz\H onyi--Weiner Lemma, {\cite[Corollary 2.4]{WSzT}, \cite[Appendix, Result 6]{Hetamas}}]
\label{SzW}
Let $f$ and $g$ be two-variable polynomials in $\GF(q)[X,Y]$. Let $d=\deg f$ and
suppose that the coefficient of $X^d$ in $f$ is non-zero.
For $y\in\GF(q)$, let $h_y=\deg\gcd\left(f(X,y),g(X,y)\right)$, where $\gcd$
denotes the greatest common divisor of the two polynomials in $\GF(q)[X]$.
Then for any $y_0\in\GF(q)$,
\[\sum_{y\in\GF(q)}{\left(h_y-h_{y_0}\right)^+} 
\leq (\deg{f(X,Y)}-h_{y_0})(\deg{g(X,Y)}-h_{y_0}).\]
\end{theorem}

A partial cover of $\PG(2,q)$ with $h>0$ holes is a set of lines in $\PG(2,q)$ such that the union of these lines
cover all but $h$ points.
We will also use the dual of the following result due to Blokhuis, 
Brouwer and Sz\H{o}nyi \cite{covering}.

\begin{theorem}[{\cite[Proposition 1.5]{covering}}]
\label{cover}
A partial cover of $\PG(2, q)$ with $h > 0$ holes, not all on one line if $h>2$, has
size at least $2q-1-h/2$.
\end{theorem}

\begin{lemma}
\label{main}
Let $\so$ be a set of points in $\PG(2,q)$, let $\ell$ be a $k$-secant of $\so$ with $2 \leq k\leq q$ and
let $1 \leq t \leq q-3$ be an integer. 
Suppose that through each point of $\ell \cap \so$ there pass exactly $t$ tangent lines to $\so$.
Denote by $s$ the size of $\so$ and let $s=k+q-t+\ep$.
Let $A(n)$ be the set of those points in $\ell \setminus \so$ through which there pass at most $n$ skew lines to $\so$.
Then the following hold. 
\begin{itemize}
\item
If $t=1$ and 
\begin{enumerate}
	\item $\ep < \frac{k}{2}-1$, then the $k$ tangent lines at the points of $\so \cap \ell$ 
	and the skew lines through the points of $A(2)$ belong to a pencil
	(hence $A(2)\setminus A(1)$ is empty),
	\item if $\ep < \frac{2k}{3}-2$, then the $k$ tangent lines 
	at the points of $\so \cap \ell$ 
	either belong to two pencils or they form a dual $k$-arc.
	If $k<q$, then the skew lines through the points of $A(2)$ 
	belong to the same pencils or dual $k$-arc.
\end{enumerate}
\item
If $t\geq 2$ and $k>q-\frac{q}{t}+1$, then 
\begin{enumerate}
\setcounter{enumi}{2}
	\item if $\ep < \frac{k}{t+1}-\frac{t}{2}$, then the $kt$ tangent lines at the points 
	of $\so \cap \ell$ and the skew lines through the points of $A(t+1)$ belong to $t$ pencils whose carriers 
	are not on $\ell$ (hence $A(t+1)\setminus A(t)$ is empty),
	\item if $\ep < \frac{k}{t+1}-1$ and $t\leq \sqrt{q}$, then the $kt$ tangent lines at 
	the points of $\so \cap \ell$ belong to $t+1$ pencils whose carriers are not on $\ell$. 
	If $k<q$, then the skew lines through the points of $A(t+1)$ belong to the same pencils.
\end{enumerate}
\end{itemize}
\end{lemma}

\proof
First we introduce some notation.
Let $\cH$ be any subset of points of the line $\ell$.
We define the line-set $\cL_{\cH}$ as $\cL_{\cH}=\{r\in\PG(2,q)\colon r\cap\ell\in
((\os\cap\ell)\cup A(t+1))\setminus \cH, r\cap(\os\setminus\ell)=\emptyset\}$,
that is the set of tangent lines to $\os$ at the points of $\so\cap\ell$ 
together with the set of skew lines to $\os$ through the points of $A(t+1)$, 
except those lines that intersect $\ell$ in a point of $\cH$. 
For each point $P\in \PG(2,q)\setminus \ell$ we define the \emph{$\cH$-index of $P$}, 
in notation $\ind_{\cH}(P)$, 
as the number of lines of $\cL_{\cH}$ that pass through $P$. Also, let 
$k_{\cH}=|(\ell\cap\os)\setminus \cH|$, $a_{\cH}=|A(t+1)\setminus \cH|$, and 
let $\delta_{\cH}$ be the number of skew lines through the points of
$A(t+1)\setminus \cH$. If $\cH=\emptyset$, we omit the prefix and the
subscript $\emptyset$, e.g.\ we write $\cL$ and $\ind(P)$ instead of
$\cL_{\emptyset}$ and $\ind_{\emptyset}(P)$. If $\cH=\{Q\}$ for some
$Q\in\ell$, we write $Q$ as prefix or subscript instead of $\{Q\}$. If
$P=(m)$, we write e.g.\ $\ind(m)$ instead of $\ind((m))$. Note that if
$Q\in\ell\setminus(\so\cup A(t+1))$, then $\ind_Q(P)=\ind(P)$ 
for all $P\in\PG(2,q)\setminus\ell$.

If $P\in \so \setminus \ell$, then the $\cH$-index of $P$ is 0 for any $\cH\subset\ell$.
Let $P\in\PG(2,q)\setminus (\ell \cup \so)$ be an arbitrary point 
and $Q\in\ell$. Choose the system of reference so that 
$P\in \ell_{\infty}\setminus \{Y_{\infty}\}$, $Q=Y_{\infty}$ 
and the points of $(\so \cap \ell) \cup A(t+1)\setminus \{Q\}$ are 
affine points on the line $[1:0:0]$.
Then $P=(y_0)$ for some $y_0 \in \GF(q)$. Let
$\{(0,c_1),\ldots,(0,c_{k_Q+a_Q})\}$ be the set of points of $(\so \cap \ell)
\cup A(t+1)\setminus \{Q\}$, let $D=(\ell_{\infty}\setminus \{Y_{\infty}\}) \cap \so$, $|D|=d$ and let $\cU=\so
\setminus (\ell \cup \ell_{\infty})=\{(a_1,b_1),\ldots,(a_{s-d-k},b_{s-d-k})\}$. 
Consider the R\'edei polynomials of 
$(\so \cap \ell)\cup A(t+1)\setminus \{Q\}$ and $\cU$.
Let us denote them by 
$f(X,Y)=\prod_{j=1}^{k_Q+a_Q} (X-c_j)$ 
and $g(X,Y)=\prod_{j=1}^{s-d-k}(Ya_j+X-b_j)$, respectively.
Let $\overline{D}=\ell_{\infty} \setminus (D\cup \{Y_{\infty}\})$. 
Then for any point $(y)\in \overline{D}$, 
\[h_y:=\deg\gcd\left(f(X,y),g(X,y)\right)=k_Q+a_Q-\text{ind}_Q\,(y).\]
Applying the Sz\H onyi--Weiner Lemma for 
the polynomials $f(X,Y)$ and $g(X,Y)$ we get 
\[\sum_{(y)\in \overline{D}}{\left(\text{ind}_Q\,(y_0)-\text{ind}_Q\,(y)\right)} \leq
\sum_{(y)\in \GF(q)}{\left(\text{ind}_Q\,(y_0)-\text{ind}_Q\,(y)\right)^+} \leq
 \text{ind}_Q\,(y_0)(s-d-k-k_Q-a_Q+\text{ind}_Q\,(y_0)).\]
After rearranging it we obtain
\begin{equation}
0\leq \text{ind}_Q\,(y_0)^2-\text{ind}_Q\,(y_0)(q+k+k_Q+a_Q-s)+\sum_{(y)\in \overline{D}}\text{ind}_Q\,(y).
\end{equation}
Because $\sum_{(y)\in \overline{D}}\text{ind}_Q\,(y)=k_Qt+\delta_Q$, hence
\begin{equation}
\label{le2}
0 \leq \text{ind}_Q\,(y_0)^2-\text{ind}_Q\,(y_0)(k_Q+a_Q+t-\ep)+k_Qt+\delta_Q.
\end{equation}

First we prove parts 1, 3 and 4 simultaneously. If we choose $Q$ so 
that $Q \in \ell \setminus \so$, then $k_Q=k$.
Thus the condition $\ep<\frac{k}{t+1}-1$ and the obvious 
fact $\delta_Q\leq (t+1)a_Q$ imply that 
(\ref{le2}) gives a contradiction for
$t+1 \leq \text{ind}_{Q}\,(y_0)\leq k+a_Q-\ep-1$.
We say that a point $P$ has large $Q$-index if 
$\ind_{Q}\,(P)\geq k+a_Q-\ep$ holds.

We are going to prove that each line of $\cL_Q$ contains a point with large
$Q$-index. First let $\ell'\in\cL_Q$ be a tangent to $\so$ at a point 
$T\in \ell\cap \so$. Suppose that each point of $\ell'$ has index at most $t$. 
Then
\begin{equation}
\label{eqi}
\sum_{P\in \ell' \setminus T} \text{ind}_Q\,(P) \leq tq.
\end{equation}
On the other hand, the sum on the left-hand side is at least $(k-1)t+q$, contradicting our assumption on $k$.
Similarly, if $\ell''$ is a skew line to $\os$ passing through a point $T\in
A(t+1)\setminus \{Q\}$, then the right-hand side of \eqref{eqi} remains the
same and the left-hand side is at least $kt+q$, that is a 
contradiction, too. Hence there are at least $t$ points with large $Q$-index.

Suppose that there are more than $t$ points with large $Q$-index and let
$R_1,R_2,\ldots,R_{t+1}$ be $t+1$ of them. Then 
\[(t+1)(k+a_Q-\ep) \leq 
\sum_{j=1}^{t+1} \text{ind}_Q\,(R_j) \leq \binom{t+1}{2}+tk+(t+1)a_Q.\]
This is a contradiction if $\ep < \frac{k}{t+1} - \frac{t}{2}$, which 
holds in parts 1 and 3. 
Regarding part 4, if there would be more than $t+1$ points with large
$Q$-index, then $\ep < \frac{2k}{t+2}-\frac{t+1}{2}$ yields a 
contradiction.
The condition on $k$ and $t\leq \sqrt{q}$ 
imply $\frac{k}{t+1}-1 \leq \frac{2k}{t+2}-\frac{t+1}{2}$.

If $k+|A(t+1)|<q+1$, then let $Q$ be any point of $\ell \setminus (\so \cup
A(t+1))$. Thus the lines of $\cL_Q=\cL$ are contained in $t$ pencils (or $t+1$
in part 4) whose carriers have large $Q$-index. In this case parts 1, 3 and 4
are proved. So from now on we assume $k+|A(t+1)|=q+1$. Then we set $Q\in
A(t+1)$. To finish the proof of parts 1, 3 and 4, we have to show that the
lines of $\cL\setminus\cL_Q$ also belong to these pencils. Recall that in case
of part 4, we assume $k<q$.

If $k=q$, then let $Q$ be the unique point contained in $A(t+1)$.
The $kt$ tangents at the points of $\ell \cap \so $ are contained 
in $t$ pencils having carriers with large $Q$-index.
Denote the set of these carriers by $\cP = \{G_1,G_2,\ldots,G_t\}$.
If $t=1$, then through $G_1$ there pass $q$ tangent lines, 
hence the points of $\so \setminus \ell$ are contained in the line $G_1Q$.
Thus through $Q$ there pass only two non-skew lines, $\ell$ and $G_1Q$. 
The condition $q-3\geq t=1$ implies $(q+1)-2 > 2$, 
hence $Q \notin A(2)$, a contradiction.
If $t>1$, then $\cP$ is contained in a line through $Q$ and 
$\so \setminus \ell$ is a $(q-t)$-secant of $\so$.
Again $q-3\geq t$ implies that through $Q$ there pass more than $t+1$ skew 
lines, hence $Q \notin A(t+1)$, a contradiction.

If $k<q$, then let $Q_1$ and $Q_2$ be two distinct points of $A(t+1)$.
As seen before, the lines of $\cL_{Q_i}$ are blocked by the points with large
$Q_i$-index for $i=1,2$, hence, by $\cL_{Q_1}\cup\cL_{Q_2}=\cL$, it is enough
to show that the set of points with large $Q_1$-index is the same as the set
of points with large $Q_2$-index. 
If a point has large $Q_i$-index, then its $Q_i$-index is at least 
$k+a_{Q_i}-\ep=q-\ep$, while the other points have $Q_i$-index at most $t$ 
for $i=1,2$. 
The inequality $|\text{ind}_{Q_1}\,(P) - \text{ind}_{Q_2}\,(P)|\leq 1$
obviously holds, thus it is enough to show that $q-\ep - t > 1$, 
which follows from the assumptions $\ep < \frac{k}{t+1}-1$ and $t\leq q-3$.

Finally, we prove part 2. At this part sometimes we will choose $Q$ from $\ell
\cap \so$, so from now on $k_Q$ is not necessarily equal to $k$. Let $P$ be
the point of $\PG(2,q)\setminus (\ell\cup \so)$ whose index is to be
estimated. If $k+|A(2)|<q+1$, then let $Q$ be any point of $\ell \setminus
(\so \cup A(2))$ and let $\cH=\emptyset$. If $k+|A(2)|=q+1$ and $k=q$, then
let $Q$ be the unique point contained in $A(2)$ and let $\cH=\{Q\}$, otherwise
let $Q$ be any point of $\ell $ such that $PQ$ intersects 
$\so \setminus \ell$ and let
$\cH=\emptyset$. Note that since $\so \setminus \ell$ is not empty, $Q$ can be
chosen in this way and $\ind_{\cH}\,(P)$ does not depend on the choice of
$Q$. In all cases we investigate the line-set $\cL_{\cH}$, and we have
$\text{ind}_Q\,(P)=\text{ind}_{\cH}\,(P)$. If
\begin{equation}
\label{part2}
\frac{2k_Q}{3}-2+\frac{a_Q}{3}>\ep
\end{equation}
holds, then \eqref{le2} gives a contradiction for 
$3 \leq \text{ind}_{\cH}\,(P)\leq k_Q+a_Q-2-\ep$.
In all cases the left-hand side of $(\ref{part2})$ is at 
least $2k/3-2$, hence the corresponding lines either form 
a dual arc or there is a point $G$ with $\cH$-index at least $k_Q+a_Q-1-\ep$.

In the latter case let $B=(\ell \setminus (\so \cup A(2))) \cup (\so \setminus \ell) \cup G$
and denote by $h$ the number of lines of $\PG(2,q)$ not blocked by $B$.
It is easy to see that, apart from $\ell$, $B$ blocks all but at 
most $(k+2|A(2)|)-(k_Q+a_Q-1-\ep)$ lines of $\PG(2,q)$.
If $k+|A(2)|<q+1$, then $B$ blocks $\ell$ and $k_Q+a_Q=k+|A(2)|$, 
hence $h\leq|A(2)|+1+\ep$.
If $k+|A(2)|=q+1$, then $B$ does not block $\ell$ and $k_Q+a_Q=q$, 
thus $h \leq |A(2)|+3+\ep$.

Suppose to the contrary that these $h$ lines do not pass through one point. 
Then from Theorem \ref{cover} we have
$|B|\geq 2q-1-h/2$ or, equivalently,
	\[q+1-(k+|A(2)|)+(q-1+\ep)+1 \geq 2q-1-h/2.\]
Rearranging it we obtain $\ep \geq k+|A(2)|-2-h/2$. 
If $k+|A(2)|<q+1$, then this would imply $\ep \geq 2k/3-5/3+|A(2)|/3$. 
If $k+|A(2)|=q+1$, then $\ep \geq (q+k)/3-2$ would follow.
Both cases yield a contradiction because of our assumption on $\ep$.
Hence the corresponding lines can be blocked by $G$ and one more point.
\Qed

\begin{corollary}
\label{so}
Let $\so _1$ be a semioval in $\PG(2,q)$ and let $\ell $ be a $k$-secant 
of $\so _1$. If $|\so _1|< q+\frac{3k}{2}-2$, then the $k$ tangent lines 
at the points of $\so _1\cap \ell$ belong to a pencil.
If $|\so _1|< q + \frac{5k}{3}-3$, then the $k$ tangent lines at the points 
of $\so _1\cap \ell$ either belong to two pencils or they form a dual $k$-arc.
\end{corollary}

If $k=q-1$, then we get a stronger result than the previous characterization 
of Kiss \cite[Corollary 3.1]{kgy}. 

\begin{corollary}
\label{cor2}
Let $\so _1$ be a semioval in $\mathrm{PG}(2,q)$. If $\so _1$ has a 
$(q-1)$-secant $\ell$ and $|\so _1|<\frac{5q}{2}-\frac{7}{2}$ holds, 
then $\so _1$ is contained in a vertexless triangle and it has two $(q-1)$-secants. 
\end{corollary}
\proof
Let $\ell \setminus \so_1 = \{A,B\}$. It follows from Corollary \ref{so} that the tangents at the points of
$\so_1 \cap \ell$ are contained in a pencil with carrier $C$. Thus $\so_1$ is contained in the sides of the triangle $ABC$.
Suppose to the contrary that each of $AC$ and $BC$ intersects $\so_1$ in less than $q-1$ points. Then there exist $P$, $Q$
such that $P\in AC \setminus (\so_1 \cup \{A,C\})$ and $Q \in BC \setminus (\so_1 \cup \{B,C\})$. The point $E:=PQ \cap AB$ is in $\so_1$ and 
$PQ$ is a tangent to $\so_1$ at $E$. This is a contradiction since $C\notin PQ$. 
\Qed

Since $t<q$ implies $k\leq q+1-t$, the assumption $q-\frac{q}{t}+1<k$ 
in Lemma \ref{main} can hold only if $t<\sqrt{q}$.

\begin{corollary}
\label{st}
Let $\st$ be a $t$-semiarc in $\PG(2,q)$, $q\geq 7$, with $1<t<\sqrt{q}$. 
Suppose that $\st$ has a $k$-secant $\ell$ and $q-\frac{q}{t}+1<k$ holds.
If $|\st|< (q-t+k) + \frac{k}{t+1}-1$, then the $kt$ tangent lines at the 
points of $\st \cap \ell$ belong to $t+1$ pencils. If 
$|\st|< (q-t+k) + \frac{k}{t+1}-\frac{t}{2}$, then the $kt$ tangent lines 
at the points of $\st \cap \ell$ belong to $t$ pencils.
\end{corollary}

\begin{remark}
Theorem \ref{Dir1} follows from Lemma \ref{main} with $t=1$ and $\ep=0$. 
To see this, let $\so=\cU \cup (\ell_{\infty} \setminus D_{\cU})$. Then through each point of $\ell_{\infty} \cap \so$, 
there passes a unique tangent to $\so$. According to Lemma \ref{main}, these tangent lines are contained in a pencil, whose 
carrier can be added to $\cU$.
\end{remark}

\begin{example}
It follows from Theorem \ref{cover} that a cover of the complement of a 
conic in $\PG(2,q)$, $q$ odd, by external lines,
contains at least $3(q-1)/2$ lines, see \cite[Proposition 1.6]{covering}. 
Blokhuis et al.\ also remark that this bound can be reached for $q=3,5,7,11$ and 
there is no other example of this size for $q<25$, $q$ odd.
Now, let $\ell$ be a tangent to a conic $\cC$ at the point $P\in \cC$ and let 
$\cU$ be a set of $3(q-1)/2$ interior points of the conic such that
these points block each non-tangent line. From the dual of the above result 
we know that such set of interior points exists at least when
$q=3,5,7,11$. Let $\so=\cU \cup \ell \setminus \{P\}$. Then the tangents 
to $\so$ at the points of $\ell \cap \so$ obviously do not pass through one 
point and this shows that part 1 of Lemma \ref{main} is sharp if $k=q$ 
and $q=5,7,11$. 
\end{example}

\begin{example}[{\cite[Theorem 3.2]{smallsemi}}]
\label{ex1}
In $\PG(2,8)$ there exists a semioval $\so_1$ of size 15, contained in a triangle
without two of its vertices.
The side opposite to the one vertex contained in $\so_1$ is a 6-secant and the
other two sides are 5-secants.
The tangents at the points of the 6-secant do not pass through one point.
Hence Corollary \ref{so} is sharp at least for $q=8$.
\end{example}

In the following we give some examples for small $t$-semiarcs with long secants in 
the cases $t=1,\,2,\,3$ such that the tangents at the points of a long secant do not belong to 
$t$ pencils. These assertions can be easily proved using Menelaus' Theorem.
Denote by $\GF(q)^+$ and $\GF(q)^\times$ the additive and multiplicative 
groups of the field $\GF(q)$, $q=p^h$, $p$ prime, respectively, and 
by $A \oplus B$ the direct sum of the groups $A$ and $B$.

\begin{example}[{\cite[p.\ 104]{smallsemi}}]
Consider $\GF(q)$, $q$ square, as the quadratic extension of $\GF(\sqrt{q})$ 
by $i$. Then the pointset 
$\so_1= [1:0:0]\cup [1:0:1]\cup [0:0:1] \setminus \{Y_{\infty},\,(0:s:1),\,(1:si:1),\,(1:s+si:0) \colon s \in \GF(\sqrt{q})\}$ 
is a semioval in $\PG(2,q)$ with three $(q-\sqrt{q})$-secants.
\end{example}

\begin{example}
Let $\GF(q)^+=A\oplus B$, where $A$ and $B$ are proper subgroups 
of $\GF(q)^+$ and let $X=A \cup B$.
The pointset $\so_2= \{(0:s:1),\,(1:s:1),\,(1:s:0) \colon s \in \GF(q) \setminus
X\}$, is a 2-semiarc in $\PG(2,q)$ with 
three $(q+1-|A|-|B|)$-secants.
Note that $2\sqrt{q}\leq |A|+|B| \leq q/p+p$.

Similarly, let $\GF(q)^\times=A\oplus B$ and $X=A \cup B$, where $A$ and $B$ 
are proper subgroups of $\GF(q)^\times$. 
The pointset $\so_3= \{(0:s:1),\,(s:0:1),\,(1:-s:0) \colon s \in \GF(q) \setminus (X \cup \{0\})\}$, is a 3-semiarc in $\PG(2,q)$ with three $(q-|A|-|B|)$-secants. Note that $2\sqrt{q}\leq |A|+|B| \leq (q+3)/2$.
\end{example}

\section{Semiarcs and blocking sets}
\label{sec:blsets}

First we associate a blocking set to each semiarc. 

\begin{lemma}
\label{makebs}
Let $\Pi_q$ be a projective plane of order $q$, let 
$k\leq q$ and $1 \leq t \leq q-3$ be integers. 
Let $\so$ be a set of $k+q-t+\ep$ points in $\Pi_q$ such that the line 
$\ell $ is a $k$-secant of $\so$. Let $A(n)$ be the set of those points in 
$\ell \setminus \so$ through which there pass at most $n$ skew lines to $\so$.
Suppose that through each of the $k$ points of $\ell \cap \so$ there pass 
exactly $t$ tangent lines to $\so$, and also suppose that these $kt$ tangent 
lines and the skew lines through the points of $A(n)$ belong to $n$ pencils. 
Let $\cP$ be the set of carriers of these pencils and assume that
$\cP \cap \ell =\emptyset$. Define the pointset $\cB_n(\so,\ell)$ in the 
following way:
\[\cB_n(\so,\ell):=\left(\ell \setminus (A(n)\cup \so)\right) \cup (\so \setminus \ell)\cup \cP.\]
Then $\cB_n(\so,\ell)$ has size $2q+1+\ep+n-t-k-|A(n)|$.
If $\ell \cap \cB_n(\so,\ell)=\emptyset$ (that is $\ell \subseteq A(n)\cup \so$), then $\cB_n(\so,\ell)$ is an affine 
blocking set in the affine plane $\Pi_q \setminus \ell$; otherwise $\cB_n(\so,\ell)$ is a blocking set in $\Pi_q$. 
In the latter case the points of $\ell \cap \cB_n(\so,\ell)$ are essential points.
\end{lemma}
\proof
Let $\ell'\neq \ell$ be any line in $\Pi_q$ and let $E$ be the point  $\ell \cap \ell'$.
If $\ell'$ meets $\left(\ell \setminus (A(n)\cup \so)\right) \cup (\so \setminus \ell)$, then $\ell'$ is blocked by $\cB_n(\so,\ell)$.
Otherwise $\ell'$ is a tangent to $\so$ at a point of $\ell\cap \so$ or $\ell'$ is a skew line to $\so$ that intersects $A(n)$.
In both cases $\ell'$ is blocked by $\cP$, hence it is also blocked by $\cB_n(\so,\ell)$.

If $\ell \subseteq A(n)\cup \so$, then $\cB_n(\so,\ell)$ is an affine 
blocking set in the affine plane $\Pi_q \setminus \ell$. Otherwise $\ell$ is 
blocked by $\ell \setminus (A(n)\cup \so)$ and hence $\cB_n(\so,\ell)$ is a 
blocking set in $\Pi_q$. In the latter case through each point of 
$\ell \cap \cB_n(\so,\ell)$ there are at least $n+1$ skew lines to $\so$ and 
hence through each of them there is at least one tangent to $\cB_n(\so,\ell)$.
\Qed

In $\PG(2,q)$ we will combine Lemma \ref{makebs} with Lemma 
\ref{main} in the cases $n=t$ or $n=t+1$. In this way we will get small 
blocking sets starting from small semiarcs having a long secant. 
The pointset $\cB_n(\so, \ell)$ is an affine blocking set if and only if
$k+|A(n)|=q+1$, and in this case $|\cB_n(\so, \ell)|=q+\ep+n-t$. An affine
blocking set in $\AG(2,q)$ has at least $2q-1$ 
points (see \cite{BS} or \cite{jamison}; also follows from 
Theorem \ref{cover}). Hence if we consider $\PG(2,q)$, then $\ep < q-n+t-1$ 
implies that $\cB_n(\so, \ell)$ is not an affine blocking set.
This condition will always hold for $n=t$ or $n=t+1$.

\begin{example}
\label{projtriangle}
If $\so_1$ is the semioval described in Example \ref{projtriangle0} 
and $\ell$ is one of the $(q-1)/2$-secants of 
$\so_1$, then $\so_1$ and $\ell$ satisfy the conditions of Lemma \ref{makebs} with $n=1$ and the obtained blocking set
$\cB_1(\so_1,\ell)$ is a minimal blocking set called the projective triangle (see e.g.\ \cite[Lemma 13.6]{jwph1}).
\end{example}

\begin{lemma}
\label{small1}
Let $\st$ be a $t$-semiarc in $\PG(2,q)$, $q=p^h$, $p$ prime, 
with $t\leq \sqrt{2q/3}$. Let $\ell$ be a $k$-secant of $\st$ and suppose 
that $\st$ and $\ell$ satisfy the conditions of Lemma \ref{makebs} with $n=t$.
If $p=2$ and $\ep<k-\frac{4}{5}(q-1)$, or $p$ is odd and $\ep < k-\frac{1}{2}(q-1)$, then $|A(t)|\geq t$.
\end{lemma}
\proof
In both cases we have $|\cB_{t}(\st,\ell)|=2q+1+\ep-k-|A(t)|<3(q+1)/2$, 
hence $\cB_{t}(\st,\ell)$ is a small blocking set.
Let $\cB$ be the unique (cf.\ Theorem \ref{modp}) minimal blocking set contained in it and let $e$ 
be the exponent of $\cB$ (cf.\ Theorem \ref{modpe}). Note that if $\ep<k-\frac{4}{5}(q-1)$, 
then $p^e\geq 8$ follows from Theorem \ref{modpe}.
Also $p^e\geq 3$ holds when $p$ is odd. 

The points of $\ell\cap\cB_t(\st,\ell)$ are essential points 
of $\cB_t(\st,\ell)$ hence $\ell \cap  \cB_t(\st,\ell)= \ell \cap \cB$.
The size of $\cB\cap (\st \setminus \ell)$ is at least $q-t$; 
let $\cU$ be $q-t$ points from this pointset.
Consider $\ell$ as the line at infinity and apply Lemma \ref{small0} 
with $E=A(t)$, $F=\ell\setminus \st$, $z=p^e$ and with $\cP$ defined as in Lemma \ref{makebs}.
Note that $t\leq \sqrt{2q/3} \leq \sqrt{q(z-1)/z}$.
Through each point of $\cU$ there pass $t$ tangents to $\st$. 
These lines are also tangents to $\cU$ and they have direction in $F$.
If $\ell'$ is one of these tangents, then $\cB\cap \ell' \equiv 1 \pmod z$ 
thus if $\ell'$ has direction in $F\setminus E$, then
$(\cP \cup \cU) \cap \ell' \equiv 0$ (mod $z$).
Hence the two required properties of Lemma \ref{small0} hold, thus $|A(t)|\geq t$.
\Qed

Semiarcs with two long secants were investigated by Csajb\'{o}k. He proved the following.

\begin{lemma}[{\cite[Theorem 13]{csb3}}]
\label{2hosszualtalaban}
Let $\Pi_q$ be a projective plane of order $q$, $1< t <q$ an integer and $\st$ be a $t$-semiarc in $\Pi_q$.
Suppose that there exist two lines
$\ell_1$ and $\ell_2$ such that $|\ell_1 \setminus (\st\cup \ell_2)|=n$ and $|\ell_2 \setminus (\st\cup \ell_1)|=m$.
If $\ell_1\cap \ell_2 \notin \st$, then $n=m=t$ or $q \leq \min\{n,m\}+2nm/(t-1)$.
\end{lemma}

The complete characterization of $t$-semiarcs in $\PG(2,q)$ with two $(q-t)$-secants whose common point is not in the semiarc was also given in \cite{csb3}. Here we cite just a particular case.

\begin{theorem}[{\cite[Theorem 22]{csb3}}]
\label{2hosszu}
Let $\st$ be a $t$-semiarc in $\PG(2,q)$, $q=p^h$, $p$ prime, with two $(q-t)$-secants such that the point of intersection of 
these secants is not contained in $\st$, and let $t\leq q-2$. Then the following hold.
	\begin{enumerate}
		\item If $\gcd(q,t)=1$, then $\st$ is contained in a vertexless triangle.
		\item If $\gcd(q,t)=1$ and $\gcd(q-1,t-1)=1$, then $\st$ consists of the symmetric difference of two lines with $t$ further points removed from each line.
		\item If $\gcd(q-1,t)=1$, then $\st$ is contained either in a vertexless triangle, or in the union of three concurrent lines with their common point removed.
	\end{enumerate}
\end{theorem}

Now we are ready to prove our main characterization theorems for small semiarcs with a long secant. We distinguish two cases, as the results on blocking sets in $\PG(2,q)$ are stronger if $q$ is a prime.

\begin{theorem}
\label{thmprimeorder}
Let $\st$ be a $t$-semiarc in $\PG(2,p)$, $p$ prime, and let $\ell$ be a $k$-secant of $\st$.
\begin{enumerate}
	\item If $t=1,\,p\geq 5$ and $k \geq \frac{p-1}{2}$, then 
	\begin{itemize}
	\item
	$\so_1$ is contained in a vertexless triangle and has two $(p-1)$-secants, or 
	\item
	$\so_1$
	is projectively equivalent to Example \ref{projtriangle0}, or 
	\item
	$|\so_1|\geq \min\{\frac{3k}{2}+p-2,\,2k+\frac{p+1}{2}\}$.
	\end{itemize}
	\item If $t=2,\, p\geq 7$ and $k \geq \frac{p+3}{2}$, then 
	\begin{itemize}
	\item
	$\so_2$ consists of the symmetric difference of two lines with two further points removed from each line, or
	\item
	$|\so_2|\geq \min\{\frac{4k}{3}+p-3,\,2k + \frac{p-1}{2}\}$.
	\end{itemize}
	\item If $3 \leq t < \sqrt{p},\, p \geq 23$ and $k > p-\frac{p}{t}+1$, then 
	\begin{itemize}
	\item
	$\st$ is contained in a vertexless triangle and has two
	$(p-t)$-secants, or 
	\item
	$|\st|\geq k\frac{t+2}{t+1}+p-t-1$.
\end{itemize}
\end{enumerate}
\end{theorem}

{\bf Proof of part 1.}
Assume that $|\so_1|< \min\{\frac{3k}{2}+p-2,\,2k+\frac{p+1}{2}\}$.
If $|\so_1|=k+p-1+\ep$, then we have $\ep < \min\{ \frac{k}{2}-1,\,
k-\frac{p-3}{2}\}$, hence Lemma \ref{main} implies that the tangents at the
points of $\ell\cap \so_1$ and the skew lines through the points of $A(1)$ are contained in
a pencil with carrier $P$. Construct the small blocking set 
$\cB_1(\so_1,\ell)$ as in Lemma \ref{makebs} with $n=1$.
The size of $\cB_1(\so_1,\ell)$ is $2p+1+\ep-k-|A(1)| < 3(p+1)/2+1$, thus 
Theorem \ref{smallbs} implies that $\cB_1(\so_1,\ell)$ either 
contains a line, or it is a minimal blocking set of size $3(p+1)/2$, each of 
its points has exactly $(p-1)/2$ tangents, and $\ep=k-\frac{p-1}{2}$.

In the first case, let $\ell_1$ be the line contained in $\cB_1(\so_1,\ell)$.
Since no $p$ points of $\so_1$ can be collinear (by Theorem \ref{q+1-t}), we have that $\ell_1$ 
is a $(p-1)$-secant of $\so_1$. The assertion now follows from Corollary 
\ref{cor2}. In the latter case, since the point $P\in \cB_1(\so_1,\ell)$ 
has at least $k$ tangents, we have $k=(p-1)/2$ and hence $\ep=0$. It follows 
from Corollary \ref{cort1} that $\so_1$ is projectively equivalent to the 
projective triangle.

{\bf Proof of part 2.}
Assume that $|\so_2|< \min\{\frac{4k}{3}+p-3,\,2k+\frac{p-1}{2}\}$.
If $|\so_2|=k+p-2+\ep$, then we have $\ep < \min\{ \frac{k}{3}-1,\,
k-\frac{p-3}{2}\}$, hence Lemma \ref{main} implies that the tangents at the
points of $\ell\cap \so_2$ and the skew lines through the points of $A(2)$ are contained in
two pencils with carriers $P_1$ and $P_2$.
Construct the blocking set $\cB_2(\so_2,\ell)$ as in Lemma \ref{makebs}.
Theorem \ref{smallbs} implies that $\cB_2(\so_2,\ell)$ either 
contains a line or it is a minimal blocking set of size
$3(p+1)/2$ and each of its points has exactly $(p-1)/2$ tangents.

Suppose that $\cB_2(\so_2,\ell)$ contains a line $\ell_1$.
Since $\so_2$ cannot have more than $p-2$ collinear points, we have 
that $\ell_1$ is a $(p-2)$-secant of $\so_2$.
Similarly we can construct $\cB_2(\so_2,\ell_1)$ and get that there is a line
$\ell_2\neq \ell_1$ and $\ell_2 \cap \ell_1 \notin \so_2$, which is also a
$(p-2)$-secant, or $\cB_2(\so_2,\ell_1)$ is a minimal blocking set of size 
$3(p+1)/2$.

In the first case Theorem \ref{2hosszu} implies that $\so_2$ consists of the
symmetric difference of two lines with two further points removed from each 
line. If this is not the case, then $\cB_2(\so_2,\ell)$ or 
$\cB_2(\so_2,\ell_1)$ is a minimal blocking set of size $3(p+1)/2$. We may 
suppose that $\cB_2(\so_2,\ell)$ is such a blocking set, hence both $P_1$ 
and $P_2$ have exactly $(p-1)/2$ tangent lines. But this is a contradiction 
since these two points together have at least $2k$ tangents, which is 
greater than $p-1$.

{\bf Proof of part 3.}
Assume that $|\st|< k\frac{t+2}{t+1}+p-t-1$.
Then $|\st|=k+p-t+\ep$, where $\ep < \frac{k}{t+1}-1 < k-\frac{p+1}{2}$, 
hence Lemma \ref{main} implies that the tangents at the points of $\ell\cap
\st$ are contained in $t+1$ pencils. Construct the small blocking set
$\cB_{t+1}(\st,\ell)$ as in Lemma \ref{makebs}. Theorem \ref{smallbs} implies 
that $\cB_{t+1}(\st,\ell)$ contains a line $\ell_1$. Note that 
$\ell_1 \cap \ell \notin \st$. Since $\st$ cannot have more than $p-t$
collinear points we have that $\ell_1$ is a $(p-t)$-secant or a 
$(p-t-1)$-secant of $\st$.

If $\ell$ were a $(p-t)$-secant or a $(p-t-1)$-secant, then Lemma 
\ref{2hosszualtalaban} would imply that both $\ell$ and $\ell_1$ are 
$(p-t)$-secants. Otherwise $|\ell\cap \st|<|\ell_1\cap \st|$ and hence the 
conditions hold also with $\ell_1$ instead of $\ell$. Constructing 
$\cB_{t+1}(\st,\ell_1)$ we get that there is a line $\ell_2\neq \ell_1$ such that 
$\ell_2 \cap \ell_1 \notin \st$ and $\ell_2$ is either a $(p-t)$-secant or a 
$(p-t-1)$-secant. Again Lemma \ref{2hosszualtalaban} implies that both 
$\ell_1$ and $\ell_2$ are $(p-t)$-secants. Since $\gcd(t,p)=1$, Theorem 
\ref{2hosszu} implies that $\st$ is contained in a vertexless triangle.
\Qed

If the projective plane $\Pi_q$ contains a Baer subplane, then there 
exist $t$-semiarcs of size $(q-\sqrt{q}-t)+(q-t)$ with a $(q-\sqrt{q}-t)$-secant, 
see Example \ref{subplane}. The first part of the following theorem states that if a
line $\ell$ intersects a $t$-semiarc $\st$ in $\PG(2,q)$, $q$ square, in at
least $k\geq q-\sqrt{q}-t$ points, $t$ is not too large and the size of $\st$ is close to $k+q-t$,
then either $\st$ is the semiarc described in Example \ref{subplane} 
or $\st$ has two $(q-t)$-secants.

\begin{theorem}
\label{thmsquareorder}
Let $\st$ be a $t$-semiarc in $\PG(2,q)$, $q=p^h$, $h>1$ if $p$ is an odd 
prime and $h\geq 6$ if $p=2$. Suppose that $\st$ has a $k$-secant $\ell$ with 
\[ k \geq \left\{\begin{array}{ll} 
q-\sqrt{q}-t & \textrm{ if $h$ is even,} \\
q-c_pq^{2/3}-t & \textrm{ if $h$ is odd,} \end{array} \right.
\] 
where $c_p=2^{-1/3}$ for $p=2,3$ and $c_p=1$ for $p>3$ (cf.\ Theorem \ref{smallbs}).
Then the following hold.

\begin{enumerate}
\item
In case of $h=2d$ and $t<(\sqrt{5}-1)(\sqrt{q}-1)/2$,
\begin{itemize}
\item
if $|\st|< 2k+\sqrt{q}$, then $\st$ has two $(q-t)$-secants whose point of intersection
is not in $\st$;
\item
if $|\st|=2k+\sqrt{q}$ and $q>9$, then either $\st$ has two $(q-t)$-secants whose point of intersection
is not in $\st$, or $\st$ is as in Example \ref{subplane}. 
\end{itemize}
\item
If $h=2d+1$, $|\st|< 2k+c_pq^{2/3}$ and $t<q^{1/3}-3/2$ (or $t<(2q)^{1/3}-2$ when 
$p=2,3$), then $\st$ has two $(q-t)$-secants whose point of intersection is not in $\st$.
\end{enumerate}
\end{theorem}
\proof
To apply Lemma \ref{main}, we need $k > q-\frac{q}{t}+1$; furthermore, $\ep < k/2-1$ for $t=1$ and $\ep <k/(t+1)-t/2$ for $t\geq 2$. 
Let us consider the condition on $k$; that on $\ep$ we treat later. 
If $q$ is a square, then $k\geq q-\sqrt{q}-t>q-\frac{q}{t}+1$ holds if $t<\Phi(\sqrt{q}-1)$, where $\Phi=\frac{\sqrt{5}-1}{2}\approx 0.618034$.
Note that if $t<\Phi(\sqrt{q}-1)$, then Theorem \ref{q+1-t} implies that $\st$ cannot have more than $q-t$ collinear points. 
If $q$ is not a square, then $t<q^{1/3}-3/2$ (or $t<(2q)^{1/3}-2$ when $p=2,3$) and  $k\geq q-c_p q^{2/3}-t$ imply $k>q-\frac{q}{t}+1$. 

Next we define $b(q)$ as follows. 
\[ b(q):= \left\{\begin{array}{ll} 
\sqrt{q} & \textrm{ if $h$ is even,} \\
c_pq^{2/3} & \textrm{ if $h$ is odd.} \end{array} \right.
\] 

For $|\st |< 2k+b(q)$, we prove the $h$ even and $h$ odd cases of the theorem simultaneously. 
Let us verify the condition of Lemma \ref{main} on $\ep$. 
If $|\st |=k+q-t+\ep$, then $|\st | < 2k+b(q)$ implies $\ep < k-q+b(q)+t$.

If $t=1$, then the upper bounds on $t$ imply $q\geq 9$ for $h=2d$, 
and $q\geq 27$ for $h=2d+1$. From these lower bounds on $q$ and from 
$k\leq q-1$ it follows that $k/2 \leq (q-1)/2 \leq q-b(q)-2$, thus $\ep < k-q+b(q)+1\leq\frac{k}{2}-1$, so
the conditions of Lemma \ref{main} hold if $t=1$.

Now suppose that $t\geq 2$. As $\ep \leq k-q+b(q)+t\leq\frac{k}{2}-1$, it is enough to prove $k-q+b(q)+t<\frac{k}{t+1}-\frac{t}{2}$.
After rearranging we get that this is equivalent to
\[k<(q-t)+\left(\frac{q-b(q)}{t}-\frac{t}{2}- b(q)-\frac{3}{2}\right),\]
thus it is enough to see (as $k\leq q-t$ holds automatically) that
\[\frac{q-b(q)}{t}-\frac{t}{2} - b(q)-\frac{3}{2}>0.\]
As a function of $t$ the left hand side decreases monotonically. It is 
positive when $t$ is maximal (under the respective assumptions), hence the conditions of Lemma \ref{main} are 
satisfied.

Construct the blocking set $\cB_t(\st,\ell)$ as in Lemma \ref{makebs}.
The conditions in Lemma \ref{small1} hold, hence the size of $A(t)$ is at 
least $t$. The size of $\cB_t(\st,\ell)$ is $2q+1+\ep-k-|A(t)|< q+b(q)+1$, thus
Theorem \ref{smallbs} implies that $\cB_t(\st,\ell)$ contains a line $\ell_1$.
Since $\st$ cannot have more than $q-t$ collinear points, we get that 
$\ell_1$ is a $(q-t)$-secant of $\st$. We have 
$|\ell\cap \st|\leq|\ell_1\cap \st|$ and hence the conditions in Lemmas 
\ref{main} and \ref{small1} hold also with $\ell_1$ instead of $\ell$.
Constructing $\cB_t(\st,\ell_1)$ we get that there exists another 
$(q-t)$-secant of $\st$, having no common point with $\ell_1 \cap \st$.

Now consider the case $h=2d$, $|\st|= 2k+\sqrt{q}$ and suppose that $\st$ 
does not have two $(q-t)$-secants. We can repeat the above arguing and get 
that $\cB_t(\st,\ell)$ or $\cB_t(\st,\ell_1)$ is a Baer subplane because of 
Theorem \ref{smallbs}. (Here, to assure $\ep<k/2-1$ in case of $t=1$, we use $q>9$.) We may suppose that $\cB_t(\st,\ell)$ is a Baer 
subplane and hence $\ep = k-q+\sqrt{q}+t$ and $|A(t)|=t$.
The size of $\ell \cap \cB_t(\st,\ell)$ is either 1 or $\sqrt{q}+1$.
In the latter case $k=q-\sqrt{q}-t$ and we get Example \ref{subplane}.
In the first case $k=q-t$. We show that this cannot occur. 
Denote by $R$ the common point of $\ell$ and $\cB_t(\st,\ell)$ and
let $P$ be any point of $\cB_t(\st,\ell) \setminus (\ell \cup \st)$.
Among the lines of the Baer subplane $\cB_t(\st,\ell)$ there are $\sqrt{q}+1$ lines incident with $P$, one of them is $PR$, which meets $\st$ 
in at least $\sqrt{q}-t > 1$ points. 
Each of the other $\sqrt{q}$ lines of the subplane intersects $\st$ in at least $\sqrt{q}+1-t> 1$ points,
thus these lines cannot be tangents to $\st$. 
But the pencil of lines through $P$ contains $k=q-t$ tangents to $\st$, one at each point of $\ell \cap \st$, too.
Thus the total number of lines through $P$ is at least $1+\sqrt{q}+q-t>q+1$, this is a contradiction. 
\Qed

\begin{flushleft}
Bence Csajb\'{o}k \\
Department of Mathematics, Informatics and Economics\\
University of Basilicata \\
Campus Macchia Romana, via dell'Ateneo Lucano \\
I-85100 Potenza, Italy

\vspace{2mm}

MTA--ELTE Geometric and Algebraic Combinatorics Research Group\\
1117 Budapest, P\'{a}zm\'{a}ny P\'eter s\'et\'any 1/C, Hungary \\
e-mail: {\sf csajbokb@cs.elte.hu}
\end{flushleft}

\begin{flushleft}
Tam\'{a}s H\'{e}ger  \\
Department of Computer Science\\
MTA--ELTE Geometric and Algebraic Combinatorics Research Group\\
E\"{o}tv\"{o}s Lor\'{a}nd University \\
1117 Budapest, P\'{a}zm\'{a}ny P\'eter s\'et\'any 1/C, Hungary \\
e-mail: {\sf heger@cs.elte.hu}
\end{flushleft}

\begin{flushleft}
Gy\"{o}rgy Kiss \\
Department of Geometry\\
MTA--ELTE Geometric and Algebraic Combinatorics Research Group \\
E\"{o}tv\"{o}s Lor\'{a}nd University \\
1117 Budapest, P\'{a}zm\'{a}ny P\'eter s\'et\'any 1/C, Hungary \\
e-mail: {\sf kissgy@cs.elte.hu}
\end{flushleft}


\small

\begin{thebibliography}{20}

\setlength{\parskip}{-1mm}


\bibitem{Ball}
S.~Ball: {\it The number of directions determined by a function over a finite field}, J. Combin. Theory Ser. A 104 (2003), 341--350.

\bibitem{bar}
D.~Bartoli: {\it On the structure of semiovals of small size}, to appear in J. Combin. Des., DOI: 10.1002/jcd.21383.

\bibitem{bat} 
L.~M.~Batten: {\it Determining sets},
Australas. J. Combin. 22 (2000), 167--176.

\bibitem{Blokhuis}
A.~Blokhuis: {\it On the size of a blocking set in $\PG(2,p)$}, Combinatorica 14 (1994), 111--114.

\bibitem{seminuclear}
A.~Blokhuis: {\it Characterization of seminuclear sets in a finite projective plane},
J. Geom. 40 (1991), 15--19.

\bibitem{BBSSZ}
A.~Blokhuis, S.~Ball, A.~E.~Brouwer, L.~Storme and T.~Sz\H{o}nyi:
{\it On the number of slopes of the graph of a function defined on a finite field}, J. Combin. Theory Ser. A 86 (1999), 187--196.

\bibitem{bss}
A. Blokhuis, A. E. Brouwer and T. Sz\H onyi: {\it The number of directions 
determined by a function $f$ on a finite field}, J. Comb. Theory Ser. A 70 (1995), 349--353.

\bibitem{covering}
A.~Blokhuis, A.~E.~Brouwer and T.~Sz\H{o}nyi: {\it Covering all points except one}, J. Algebraic Combin. 32 (2010) 59--66. 

\bibitem{smallb}
A.~Blokhuis, L.~Storme and T.~Sz\H{o}nyi: {\it Lacunary polynomials, multiple blocking sets and Baer subplanes}, 
J. London Math. Soc. 60 (1999),
321--332.

\bibitem{BS}
A.~E.~Brouwer and A.~Schrijver: {\it The blocking number of an affine space},
J. Combin. Theory Ser. A 24 (1978), 251--253.

\bibitem{bruen}
A.~A.~Bruen: {\it Baer subplanes and blocking sets}, Bull. Amer. Math. Soc. 76 (1970), 342--344.


\bibitem{csb3} B.~Csajb\'{o}k: {\it Semiarcs with long secants}, Electron. J. Combin. 21 (2014), \# P1.60, 14 pages.



\bibitem{csk}
B.~Csajb\'{o}k and Gy.~Kiss: {\it Notes on semiarcs}, Mediterr. J. Math. 9 (2012), 677--692.

\bibitem{dov}
J.~M.~Dover: {\it Semiovals with large collinear subsets}, J. Geom. 69 (2000), 58--67.

\bibitem{Gacs}
A.~G\'{a}cs: {\it On a generalization of R\'{e}dei's theorem}, Combinatorica 23 (2003), 585--598.

\bibitem{p2}
A.~G\'{a}cs, L.~Lov\'{a}sz and T.~Sz\H{o}nyi: {\it Directions in $\AG(2,p^2)$}, Innov. Incidence Geom. 6/7 (2009), 189--201.

\bibitem{t-arc}
A.~G\'acs and Zs.~Weiner: {\it On $(q+t,t)$-arcs of type $(0,2,t)$}, Des. Codes Cryptogr. 29 (2003), 131--139.

\bibitem{Hetamas}
T.~H\'eger: {\it Some graph theoretic aspects of finite geometries}, PhD Thesis, E\"otv\"os Lor\'and University (2013).

\bibitem{jwph1}
J.~W.~P.~Hirschfeld: {\it Projective Geometries over Finite Fields}, 
2$^{nd}$ ed., Clarendon Press, Oxford, 1998.

\bibitem{jamison}
R.~E.~Jamison:
{\it Covering finite fields with cosets of subspaces},
J. Combin. Theory Ser. A 22 (1977), 253--266.

\bibitem{kgy}
Gy.~Kiss: {\it Small semiovals in $\PG(2,q)$},
J. Geom. 88 (2008), 110--115.

\bibitem{kgy2}
Gy.~Kiss: {\it A survey on semiovals}, Contrib. Discrete Math. 3 (2008), 81--95.

\bibitem{smallsemi}
Gy.~Kiss and J.~Ruff: {\it  Notes on small semiovals}, Ann. Univ. Sci. Budapest. E\"{o}tv\"{o}s Sect. Math. 47 (2004), 97--105.

\bibitem{KMarcs}
G.~Korchm\'{a}ros and F.~Mazzocca: {\it On $(q+t,t)$-arcs of type $(0,2,t)$ in a Desarguesian plane of order $q$}
, Math. Proc. Cambridge Philos. Soc. 108 (1990), 445--459.

\bibitem{LS}
L.~Lov\'asz and A.~Schrijver: {\it Remarks on a theorem of R\'edei}, Studia Scient. Math. Hungar. 16 (1981),
449--454.

\bibitem{Polverino}
O.~Polverino: {\it Small minimal blocking sets and complete $k$-arcs in $\PG(2, p^3)$}, Discrete Math. 208/209 (1999), 469--476.

\bibitem{linconj}
P.~Sziklai: {\it On small blocking sets and their linearity}, J. Combin. Theory, Ser. A 115 (2008), 1167--1182.

\bibitem{DirSziklai}
P.~Sziklai: {\it Subsets of $\GF(q^2)$ with $d$-th power differences}, Discrete Math. 208/209 (1999), 547--555.

\bibitem{Szonyi}
T.~Sz\H{o}nyi: {\it Blocking Sets in Desarguesian Affine and Projective Planes}, Finite Fields Appl. 3 (1997), 187--202.

\bibitem{DirSzonyi}
T.~Sz\H{o}nyi: {\it On the number of directions determined by a set of points in an affine Galois plane}, J. Combin. Theory Ser. A 74 (1996), 141--146.

\bibitem{WSzT}
T.~Sz\H{o}nyi and Zs.~Weiner: {\it Proof of a conjecture of Metsch},
J. Combin. Theory Ser. A 118 (2011), 2066--2070.

\bibitem{PV}
P.~Vandendriessche: {\it Codes of Desarguesian projective planes of even order, projective triads and $(q+t,t)$-arcs of type $(0,2,t)$},
Finite Fields Appl. 17 (2011), 521--531.

\end{thebibliography}

\normalsize
\end{document}